\documentclass{hjm1}
\usepackage{verbatim}
\usepackage{amssymb}
\def\id{\operatorname{id}}
\def\range{\operatorname{range}}
\def\rank{\operatorname{rank}}
\def\Pspan{\operatorname{span}}
\def\text#1{{\hbox{#1}}}
\def\qed{\nobreak\hfill\penalty250 \hbox{}\nobreak\hfill\qedbox}
\def\operatorname#1{{\rm#1\,}}
\def\text#1{{\hbox{#1}}}
\def\demo#1{\smallbreak\noindent\bf#1:\rm\ }
\def\enddemo{\smallbreak}
\def\qedbox{\hbox{$\rlap{$\sqcap$}\sqcup$}}
\def\cal{\mathcal}
\journalvol{76}
\journalno{1}
\journalyear{2050}

\begin{document}
\title[Algebraic curvature tensors]
{ Algebraic curvature tensors for indefinite metrics whose
skew-symmetric curvature operator has constant Jordan normal 
form}
\author[P. Gilkey and T. Zhang]
{Peter B. Gilkey and Tan Zhang}
\address [Peter B. Gilkey] 
{Mathematics Department, University of Oregon, Eugene Or 97403 USA} 
\email{gilkey@darkwing.uoregon.edu}\
\address[Tan Zhang]
{Department of Mathematics and Statistics,
Murray State University,
Murray, KY 42071-0009} 
\email{tan.zhang@murraystate.edu}
\cgw
\thanks{Research of the first author partially supported by the NSF
(USA) and the MPI(Leipzig)}
\begin{abstract}
We classify the connected pseudo-Riemannian manifolds
of signature $(p,q)$ with $q\ge5$ so that at each point of $M$
the skew-symmetric curvature
operator has constant rank $2$ and constant Jordan normal form
on the set of spacelike $2$ planes and
so that the skew-symmetric curvature operator is not nilpotent for at least
one point of $M$.
\end{abstract}
\maketitle
\newtheorem{theorem}{Theorem}[section]
\newtheorem{corollary}[theorem]{Corollary}
\newtheorem{lemma}[theorem]{Lemma}
\newtheorem{remark}[theorem]{Remark}

\section{\label{qctSA}Introduction}
Let $(M,g)$ be a
smooth connected pseudo-Riemannian manifold of signature $(p,q)$. Let
$\nabla$ be the Levi-Civita connection on $TM$. The {\it Riemann
curvature operator} and the associated {\it curvature tensor} are
defined by:
\begin{eqnarray*}
   &&{}^gR(x,y):=\nabla_x\nabla_y-\nabla_y\nabla_x-\nabla_{[x,y]}\text{ and}\\
   &&{}^gR(x,y,z,w):=g({}^gR(x,y)z,w).\end{eqnarray*}
We have the following symmetries
\begin{eqnarray*}
     &&{}^gR(x,y,z,w)={}^gR(z,w,x,y)=-{}^gR(y,x,z,w),\text{ and}\\
     &&{}^gR(x,y,z,w)+{}^gR(y,z,x,w)+{}^gR(z,x,y,w)=0.\end{eqnarray*}
Let $V$ be a vector space with an inner product $(\cdot,\cdot)$ of
signature $(p,q)$. We say that a $4$ tensor $R\in\otimes^4V^*$ is an
{\it algebraic curvature tensor} if it satisfies the symmetries given above; the associated
curvature operator is then defined by $R(x,y,z,w)=(R(x,y)z,w)$. We say that $(M,g)$ is a
{\it geometrical realization} of an algebraic curvature tensor $R$ at a point $P$ of $M$
if there exists an isometry $\Theta$ from $T_PM$ to $V$ so that
${}^gR(x,y,z,w)=R(\Theta x,\Theta y,\Theta z,\Theta w)$
for all $x,y,z,w\in T_PM$.
Every algebraic curvature tensor has such a geometrical realization. Conversely, it is
often useful to study problems in geometry by first passing to the purely algebraic setting
and then drawing geometrical conclusions from results obtained
algebraically.

Let $R$ be an algebraic curvature tensor on a vector space $V$ of
signature $(p,q)$. Let $\{e_1,e_2\}$ be an oriented basis for a
non-degenerate oriented $2$ plane $\pi\subset V$. Let 
$$R(\pi):=|(e_1,e_1)(e_2,e_2)-(e_1,e_2)^2|^{-1/2}R(e_1,e_2)$$
be the associated {\it skew-symmetric curvature operator}; $R(\pi)$
depends on the orientation of
$\pi$ but is independent of the particular oriented basis chosen. In this
paper, we will examine the geometric consequences which follow from assuming that the
skew-symmetric curvature operator has certain algebraic properties.

We say that a vector $v\in V$ is
{\it spacelike} if $(v,v)>0$, {\it timelike} if $(v,v)<0$, and {\it null} if $(v,v)=0$. We
say that a $2$ plane $\pi$ is {\it spacelike} if the induced metric on $\pi$ has signature
$(0,2)$, {\it timelike} if the induced metric has signature $(2,0)$, and {\it mixed} if
the induced metric has signature $(1,1)$. Otherwise $\pi$ is said to be {\it degenerate}.

The simplest invariant of a linear map is the rank. We say that an algebraic curvature
tensor $R$ has {\it spacelike rank $r$} if $\rank(R(\pi))=r$ for every
oriented spacelike $2$ plane $\pi$. The notions of
{\it timelike} and {\it mixed rank $r$} are defined similarly. If $\rank(R(\pi))$ is not constant on the set of oriented
 spacelike $2$ planes, then we say $R$ does not have constant spacelike rank. In \S\ref{qctSB},
we shall construct algebraic curvature tensors which have spacelike rank $2$ but which do
not have constant timelike or constant mixed rank.

The following result of Gilkey, Leahy and Sadofsky 
\cite{refGLS} and of Zhang
\cite{refZ} uses results from Adams \cite{refA} and Borel \cite{refBo}
to bound the rank:
\begin{theorem}\label{arefb}
Let $R$ be an algebraic curvature tensor on a vector space $V$ of
signature $(p,q)$ which has spacelike rank $r$.
\begin{enumerate}
\smallskip\item Let $p=0$. Let $q\ge 5$ and $q\ne 7,8$. Then $r\le 2$.
\smallskip\item Let $p=1$. Let $q=5$ or $q\ge9$. Then $r\le 2$.
\smallskip\item Let $p=2$. Let $q\ge10$. Then $r\le 4$. Furthermore, if
neither $q$ nor $q+2$ are powers of $2$, then $r\le 2$.
\end{enumerate}\end{theorem}

Let $V$ be a vector space of signature $(p,q)$ and let $\phi$ be a self-adjoint
linear map of $V$. We define a $4$ tensor and associated
endomorphism:
\begin{eqnarray*}
&&R_\phi(x,y,z,w):=(\phi y,z)(\phi x,w)-(\phi x,z)(\phi y,w)\text{ and}\\
&&R_\phi(x,y)z:=(\phi y,z)\phi x-(\phi x,z)\phi y.\end{eqnarray*}
The tensor $R_\phi$ is an algebraic curvature tensor and the set of tensors arising in this
fashion spans the space of all algebraic curvature tensors \cite{refGilkeyIvanova}. Such
tensors arise as the curvature tensors of hypersurfaces in flat space - see Lemma \ref{crefa}.

These
tensors play a crucial role in the classification of the algebraic curvature tensors of
rank $2$. We refer to
\cite{refGLS,refGZ} for the proof of the following result:

\begin{theorem}\label{arefc} Let $R$ be an algebraic curvature
tensor on a vector space $V$ of signature $(p,q)$ where $q\ge5$. Then $R$ has
spacelike rank $2$ if and only if $R=\pm R_\phi$ where $\phi$ is a self-adjoint
map of $V$ whose kernel contains no spacelike vectors.
\end{theorem}

The spectrum (i.e. the eigenvalues counted with multiplicity) is also a useful invariant
of a linear map. We say that an algebraic curvature tensor is {\it spacelike IP} if the
eigenvalues of
$R(\pi)$ are the same for any two oriented spacelike $2$ planes; the notions of
{\it timelike} and {\it mixed IP} are defined similarly. (The notation `IP' is chosen as
the fundamental classification results in dimension $4$ are due to Ivanov and Petrova
\cite{refIP} - see also related work in \cite{refRIGS}). One can use analytic
continuation to see that the notions of spacelike, timelike, and mixed IP coincide so we
shall simply say that
$R$ is IP if any of these three equivalent conditions holds. Two linear maps $T$ and
$\tilde T$ of $V$ are said to be {\it Jordan equivalent} if any of the
following three equivalent conditions are satisfied:
\begin{enumerate}
\item There exist bases ${\cal{B}}=\{e_1,...,e_m\}$ and 
$\tilde{\cal{B}}=\{\tilde e_1,...,\tilde e_m\}$ for
$V$ so that the matrix representation of $T$ with respect to the basis
${\cal{B}}$ is equal to the matrix representation of $\tilde T$ with
respect to the basis
$\tilde{\cal{B}}$.
\item There exists an isomorphism $\Theta$ of $V$ so $T=\Theta\tilde
T\Theta^{-1}$, i.e. $T$ and $\tilde T$ are {\it conjugate}.
\item The real Jordan normal forms of $T$ and $\tilde T$ are equal.
\end{enumerate}
\medbreak\noindent In the positive definite setting, the spectrum of a
skew-symmetric linear map determines the conjugacy class of the map. This
is not, however, the case in the indefinite setting. We say that $R$ is
{\it spacelike Jordan IP} if the Jordan normal form of $R(\cdot)$ is the
same for any two oriented spacelike $2$ planes; the notions of {\it timelike} and
{\it mixed Jordan IP} are defined similarly. We note that if $R$ is spacelike Jordan IP,
then
$R$ has spacelike rank $r$ for some $r$.

The Jordan form of a linear map determines the
spectrum. Thus if
$R$ is spacelike Jordan IP, then $R$ is spacelike IP and hence IP. In \S\ref{qctSB}, we
construct algebraic curvature tensors which are IP but which are not spacelike,
timelike, or mixed Jordan IP. We will also construct algebraic curvature tensors which
are spacelike Jordan IP but not timelike or mixed Jordan IP; thus these
notions are distinct.

Let $\phi$ be a linear map of a vector space $V$ of signature $(p,q)$. If
$(\phi v,\phi w)=(v,w)$ for all $v,w\in V$, then $\phi$ is an
{\it isometry}. If $(\phi v,\phi w)=-(v,w)$ for all $v,w\in V$, then $\phi$
is a {\it para-isometry}; para-isometries exist if and only if
$p=q$, i.e. if we are in the {\it balanced} setting. Suppose that $\phi$ is
self-adjoint. Then $\phi$ is an isometry if and only if $\phi^2=\id$;
$\phi$ is a para-isometry if and only if $\phi^2=-\id$.

We say that $R$ is {\it spacelike rank 2 Jordan IP} if $R$ has spacelike rank $2$ and if
$R$ is spacelike Jordan IP, the notions of {\it timelike rank 2 Jordan IP} and {\it mixed
rank 2 Jordan IP} are defined similarly. We have the following classification result for
such tensors
\cite{refGZ}:

\begin{theorem}\label{arefd} Let $R$ be an algebraic curvature tensor on a
vector space of signature $(p,q)$ where $q\ge5$. Then $R$ is a spacelike rank
$2$ Jordan IP algebraic curvature tensor if and only if exactly one of the following three
conditions is satisfied:\begin{enumerate}
\smallskip\item $R=CR_\phi$ where $\phi$ is a self-adjoint isometry, and where $C\ne0$.
\smallskip\item $R=CR_\phi$ where $\phi$ is a self-adjoint para-isometry, and where
$C\ne0$.
\smallskip\item $R=\pm R_\phi$ where $\phi$ is self-adjoint, where $\phi^2=0$, and where
$\ker\phi$ contains no spacelike vectors. 
\end{enumerate}\end{theorem}

\begin{remark}\label{arefe}\rm The map $\phi$ in Theorem \ref{arefd} is uniquely
defined up to sign; the constant $C$ in assertions (1) and (2) is uniquely determined. The
tensors $CR_\phi$ in assertions (1) and (2) where $\phi^2=\pm\id$ are also timelike and
mixed rank
$2$ Jordan IP. The tensor in assertion (3) is {\it nilpotent}; if $\phi^2=0$, then
$R_\phi(x,y)^2=0$ for all
$(x,y)$. This tensor is timelike rank $2$ Jordan IP if and only if $\ker\phi$ contains no
timelike vectors; it is not constant mixed rank.\end{remark}

 We say that a pseudo-Riemannian manifold $(M,g)$ is {\it spacelike rank $r$
Jordan IP} if ${}^gR$ is spacelike rank $r$ Jordan IP at every point of
$M$; the Jordan form is allowed to vary with the point but the rank is assumed to be
constant. The notions of {\it timelike rank $r$ Jordan IP} and {\it mixed rank $r$ Jordan
IP} are defined similarly. In
\S\ref{qctSC}, we will construct two families of spacelike, timelike, and mixed rank $2$ Jordan
IP pseudo-Riemannian manifolds. Lemma
\ref{crefb} deals with the pseudo-spheres and Lemma \ref{crefc} deals with warped products of a
manifold with constant sectional curvature with an interval $I\subset\mathbb{R}$.

The following is the main
result of this paper; it generalizes previously known results
\cite{refGilw,refGLS,refIP} from the Riemannian to the pseudo-Riemannian setting.

\begin{theorem}\label{areff} Let $(M,g)$ be a connected spacelike rank $2$
Jordan IP pseudo-Riemannian manifold of signature $(p,q)$ where $q\ge5$. Assume that
${}^gR$ is not nilpotent for at least one point $P$ of $M$. 
\begin{enumerate}
\smallskip\item For each point $P\in M$, we have ${}^gR_P=C(P)R_{\phi(P)}$ where $\phi(P)$
is self-adjoint map of $T_PM$ so that
$\phi(P)^2=\id$; ${}^gR_P$ is never
nilpotent.
\smallskip\item If $\phi=\pm\id$, then $(M,g)$ has constant sectional curvature and is
locally isometric to one of the manifolds constructed in Lemma {\rm\ref{crefb}}.
\smallskip\item If $\phi\ne\pm\id$, then $(M,g)$ is locally isometric to one
of the warped product manifolds constructed in Lemma {\rm\ref{crefc}}.
\end{enumerate}\end{theorem}

We shall prove Theorem \ref{areff} in \S\ref{qctSD}. The
classification of spacelike rank $2$ Jordan IP pseudo-Riemannian manifolds with nilpotent
algebraic curvature tensors is incomplete and in \S\ref{qctSE}, we present some preliminary
results on this case.

\section{\label{qctSB}Examples of algebraic curvature
tensors}

We begin this section with a technical observation.
\begin{lemma}\label{brefA} Let $V$ be a vector space with an inner
product of signature $(p,q)$. Let $\{v_1,...,v_k\}$ be a set of linearly
independent elements of $V$. Then there exist elements
$\{w_1,...,w_k\}$ of
$V$ so that $(v_i,w_j)=\delta_{ij}$.\end{lemma}

\demo{Proof} We use the inner product to define a linear map
$\psi:V\rightarrow V^*$ by the identity $\psi(w)(v)=(v,w)$. If
$w\ne0$, then there exists $v$ so $(w,v)\ne0$ and thus $\psi$ is
injective. Since $\dim V=\dim V^*$, $\psi$ is a linear
isomorphism. Let $\{v_1,...,v_k\}$ be a set of linearly independent elements
of $V$. We extend this set to a basis
$\{v_1,...,v_{p+q}\}$ for $V$ to assume without loss of generality
that $k=p+q$. Let $\{v^1,...,v^{p+q}\}$ be the associated dual basis for
$V^*$; this means if $v\in V$, then $v=\sum_iv^i(v)v_i$. The desired
elements $w_j$ of $V$ are then defined by $w_j=\psi^{-1}(v^j)$ since
$v^i(v_j)=\delta_{ij}$.\qed\enddemo

Let $V$ be a vector space of signature $(p,q)$. We can choose a {\it normalized}
orthonormal basis
$\cal{B}$ for $V$ of the form:
\begin{eqnarray}
&&\cal{B}:=\{e_1^-,...,e_p^-,e_1^+,...,e_q^+\}\text{ so }\nonumber\\
&&V^-:=\Pspan\{e_1^-,...,e_p^-\}\text{ and }\\
&&V^+:=\Pspan\{e_1^+,...,e_q^+\}
\label{brefAa}\end{eqnarray}
are maximal orthogonal timelike and spacelike subspaces.  Let $R=R_\phi$ where $\phi$ is a
self-adjoint linear map. Let $\pi$ be a non-degenerate $2$ plane. By Lemma \ref{brefA}:
\begin{eqnarray*}
   &&\rank(R_\phi(\pi))=2\text{ if }\ker\phi\cap\pi=\{0\}\text{ and}\\
   &&\rank(R_\phi(\pi))=0\text{ if }\ker\phi\cap\pi\ne\{0\}.\end{eqnarray*}
We use this observation to
show that the notions of constant spacelike, timelike, and mixed rank are
distinct notions.
\begin{lemma}\label{brefa} Let $V$ have signature $(p,q)$ where $p\ge2$ and
$q\ge2$.\begin{enumerate}
\smallskip\item There exists an algebraic curvature
tensor
$R$ on $V$ which has spacelike rank $2$ but which does not
have constant timelike or mixed rank.
\smallskip\item There exists an algebraic curvature tensor $R$ on
$V$ which has timelike rank $2$ but which does not have constant
spacelike or mixed rank.
\smallskip\item There exists an algebraic curvature tensor $R$ on
$V$ which has spacelike and timelike rank $2$ but which
does not have constant mixed rank.
\end{enumerate}\end{lemma}

\demo{Proof} Let $\cal{B}$ be a normalized
basis for $V$. Let $\phi$ be orthogonal projection on $V^+$;
$$\phi(e_i^-)=0\text{ and }\phi(e_i^+)=e_i^+.$$
As $\ker\phi=V^-$ is timelike, $\ker\phi$ contains no spacelike
vectors and
$\rank R_\phi(\pi)=2$ for any spacelike $2$ plane. Thus $R$ has spacelike
rank $2$. We define:
\begin{eqnarray*}
     &&\pi_1:=\Pspan\{2e_1^-+e_1^+,2e_2^-+e_2^+\},\ 
     \tilde\pi_1:=\Pspan\{e_1^-,e_2^-\},\\
     &&\pi_2:=\Pspan\{2e_1^-+e_1^+,e_2^+\},\qquad\phantom{aa.}
      \tilde\pi_2:=\Pspan\{e_1^-,e_2^+\}.\end{eqnarray*}
The planes $\pi_1$ and $\tilde\pi_1$ are timelike and the planes $\pi_2$ and $\tilde\pi_2$
are mixed. We show that $R$ does not have constant timelike or mixed rank by noting:
$$\ker\phi\cap\pi_i\ne\{0\}\text{ and }\ker\phi\cap\tilde\pi_i=\{0\}\text{ for }i=1,2.$$
Assertion (1) now
follows; we interchange the roles of $+$ and $-$ to prove assertion
(2) similarly.

To prove assertion (3), we define a self-adjoint linear map $\tilde\phi$ of $V$ by
setting:
\begin{eqnarray*}
&&\tilde\phi(e_1^-)=e_1^-+e_1^+,\quad \tilde\phi(e_i^-)=e_i^-\text{ for }i>1\\
&&\tilde\phi(e_1^+)=-e_1^--e_1^+,\ \tilde\phi(e_j^+)=e_j^+\text{ for }j>1.\end{eqnarray*}
Then $\ker\tilde\phi=\Pspan\{e_1^-+e_1^+\}$. Since $\ker\tilde\phi$ is totally isotropic,
$\ker\tilde\phi$ contains no spacelike vectors and no timelike vectors so
$R_{\tilde\phi}$ has constant spacelike and timelike rank $2$. We define
mixed $2$ planes 
$$\pi_1:=\Pspan\{e_1^-,e_1^+\}\text{ and }
  \pi_2:=\Pspan\{e_2^-,e_2^+\}.$$
Since $\ker\tilde\phi\cap\pi_1\ne\{0\}$ and
$\ker\tilde\phi\cap\pi_2=\{0\}$,
$R_{\tilde\phi}$ does not have constant mixed rank. Assertion (3)
now follows. \qed\enddemo

Let $R$ be an algebraic curvature tensor of constant spacelike rank $r$ on a
vector space of signature $(p,q)$. In Theorem
\ref{arefb}, we noted that $r=2$ in many cases if $p=0$, $p=1$, or $p=2$.
There are, however, examples where $r=4$ if $p\ge q$.

\begin{lemma}\label{brefb} Let $V$ be a vector space of signature
$(p,q)$ where $p\ge q\ge2$.
\begin{enumerate}
\smallskip\item If $p=q$, then there exists an algebraic curvature
tensor $R$ on $V$ which has constant spacelike and timelike rank $4$, and which
does not have constant mixed rank.
\smallskip\item If $p>q$, then there exists an algebraic curvature
tensor $R$ on $V$ which has constant spacelike rank $4$, and which does
not have constant timelike or mixed rank $4$.
\end{enumerate}\end{lemma}

\demo{Proof} Let $\cal{B}$ be a normalized basis for $V$. We
suppose that $p\ge q$ and define a linear self-adjoint map $\phi$ of
$V$ by setting:
\begin{eqnarray*}
   &&\phi(e_i^-)=-e_i^+\text{ for }1\le i\le q,\\
   &&\phi(e_i^-)=\phantom{-e}0\text{ for }q<i\le p,\text{ and}\\
   &&\phi(e_j^+)=\phantom{-}e_j^-\text{ for }1\le j\le q.\end{eqnarray*}
Expand any vector $v\in V$ in the form
$v=\textstyle\sum_ic_ie_i^-+\textstyle\sum_jd_je_j^+$. Then
$$(v,v)=-\textstyle\sum_ic_i^2+\textstyle\sum_jd_j^2\text{ and }
 (\phi v,\phi v)=\textstyle\sum_{i\le q} c_i^2
     -\textstyle\sum_jd_j^2.$$
Thus if $v$ is spacelike, $\phi v$ is timelike. Furthermore, if
$p=q$, then $\phi$ is a para-isometry. Let
$R:=R_{\id}+R_\phi$. Then
$$R(x,y)z=(y,z)x-(x,z)y+(\phi y,z)\phi x-(\phi x,z)\phi y.$$
If $\{x,y\}$ spans a spacelike $2$ plane $\pi$, then $\{\phi x,\phi  y\}$ spans a timelike
$2$ plane $\phi\pi$. Thus
$\{x,y,\phi x,\phi y\}$ is a linearly independent set and we may use
Lemma \ref{brefA} to see that $R$ has spacelike rank $4$. If
$p=q$, then the roles of $+$ and $-$ are symmetric so $R$ also has
timelike rank $4$. If $p>q$, then $\ker\phi$ contains a timelike
vector and the same argument used to prove Lemma \ref{brefa} shows $R$
does not have constant timelike rank. Let
$\pi_1:=\Pspan\{e_1^-,e_1^+\}$ and $\pi_2:=\Pspan\{e_1^-,e_2^+\}$ be
mixed $2$ planes. We show that $R$ does not have mixed rank $4$ by noting that:
$$\rank R_\phi(\pi_1)\le2\text{ and }\rank R_\phi(\pi_2)=4.\ \qed$$
\enddemo

We omit the proof of the following result as the proof is straightforward.

\begin{lemma}\label{brefc} Let $V$ be a vector space of signature
$(p,q)$.\begin{enumerate}
\smallskip\item If $T$ is a linear map of $V$ with $T^2=0$, $0$ is
the only eigenvalue of $T$.
\smallskip\item If $T_1$ and $T_2$ are linear maps of $V$ with
$T_1^2=T_2^2=0$, then $T_1$ is Jordan equivalent to $T_2$ if and only
if $\rank(T_1)=\rank(T_2)$.
\end{enumerate}\end{lemma}

We use Lemma \ref{brefc} to show that the concepts IP, spacelike Jordan IP, and
timelike Jordan IP are inequivalent.
\begin{lemma}\label{brefd} Let $V$ have signature $(p,q)$ where $p\ge3$ and
$q\ge3$.\begin{enumerate}
\smallskip\item There exists an algebraic
curvature tensor
$R$ on $V$ which is IP, not
spacelike Jordan IP, not timelike Jordan IP, and not mixed Jordan IP.
\smallskip\item If $p=q$, then there exists an algebraic curvature tensor
$R$ on
$V$ which is spacelike Jordan IP, timelike Jordan IP, and not mixed
Jordan IP.
\smallskip\item If $p>q$, then there exists an algebraic curvature tensor
$R$ on $V$ which is spacelike Jordan IP, not timelike Jordan IP, and not
mixed Jordan IP.\end{enumerate}\end{lemma}
\demo{Proof} Let $\cal{B}$ be a normalized basis for $V$. For $k\le q$, we define:
\begin{eqnarray*}
  &&\phi_k(e_i^-)=\phantom{-}e_i^-+e_i^+\text{ for }i\le k,\  
   \phi_k(e_i^-)=0\text{ for }i>k,\\
  &&\phi_k(e_i^+)=-e_i^--e_i^+,\text{ for }i\le k,\ 
   \phi_k(e_i^+)=0\text{ for }i>k.\end{eqnarray*}
Then $\phi_k$ is self-adjoint and $\range\phi_k$ is totally isotropic, i.e. the metric is
trivial on
$\range\phi_k$. Since $(\phi_ku,\phi_kv)=0$ for all $u,v\in V$, we have
$R_{\phi_k}(x,y)^2=0$ for any $x,y\in V$. Since
$R_{\phi_k}(\pi)$ is nilpotent, we apply Lemma \ref{brefc} to see that $0$ is the
only eigenvalue of $R_{\phi_k}$ and hence
$R$ is IP. We set $k=2$ to prove assertion (1). Let
\begin{eqnarray*}
 &&\pi_1:=\Pspan\{e_1^-,e_2^-\},\ \pi_2:=\Pspan\{e_1^+,e_2^+\},\ 
  \pi_3:=\Pspan\{e_1^-,e_2^+\},\\
 &&\tilde\pi_1:=\Pspan\{e_1^-,e_3^-\},\ \tilde\pi_2:=\Pspan\{e_1^+,e_3^+\},\ 
  \tilde\pi_3:=\Pspan\{e_1^-,e_3^+\}.\end{eqnarray*}
We have $\{\pi_1,\tilde\pi_1\}$ are timelike, $\{\pi_2,\tilde\pi_2\}$
are spacelike, and $\{\pi_3,\tilde\pi_3\}$ are mixed $2$ planes. For
$1\le i\le 3$, we apply Lemma
\ref{brefA} to prove assertion (1) by observing:
\begin{eqnarray*}
  &&\pi_i\cap\ker\phi_2=\{0\}\text{ so }\rank R_{\phi_2}(\pi_i)=2\\
  &&\tilde\pi_i\cap\ker\phi_2\ne\{0\}\text{ so }\rank R_{\phi_2}(\tilde\pi_i)<2.\end{eqnarray*}

Suppose $p\ge q$. We set $k=q$ to prove assertions (2) and (3). Since
$\ker\phi_q$ contains no spacelike vectors, $\rank R_{\phi_q}(\pi)=2$ for any
spacelike $2$ plane so by Lemma \ref{brefc}, $R_{\phi_q}$ is spacelike Jordan IP. Since
$\ker\phi_q$ contains no timelike vectors if and only if $p=q$, $R_{\phi_q}$ is
timelike Jordan IP if and only if $p=q$. We study
$\pi_1:=\Pspan\{e_1^-,e_2^+\}$ and
$\pi_2:=\Pspan\{e_1^-,e_1^+\}$ to see $R_{\phi_q}$ is never mixed Jordan IP.
\qed\enddemo

\section{\label{qctSC} IP Manifolds}
Let $V$ be a vector space of signature $(r,s)$, let $M$ be a simply
connected smooth manifold of dimension $m=r+s-1$, and let $F:M\rightarrow\mathbb{R}^{(r,s)}$
be an immersion of $M$. We assume the induced metric $g$ on $M$ is non-degenerate; $(M,g)$
is said to be a {\it non-degenerate hypersurface}. Let $\nu$ be a unit normal along $M$. If
$(\nu,\nu)=+1$, then $(M,g)$ has signature $(r,s-1)$. If $(\nu,\nu)=-1$, then $(M,g)$ has
signature $(r-1,s)$. We define the {\it second fundamental form}
$L$ and {\it shape operator} $S$ by:
$$L(x,y):=(xyF,\nu)\text{ and }(Sx,y):=L(x,y).$$
The following is well known so we omit the proof in the interests of brevity.

\begin{lemma}\label{crefa} Let $V$ be a vector space of signature $(r,s)$, let $(M,g)$
be a non-degenerate hypersurface in
$V$, let $\nu$ be a unit normal along $M$, and let
$S$ be the associated shape operator. Then
${}^gR=(\nu,\nu)R_S$.\end{lemma}

We say that a pseudo-Riemannian manifold $(M,g)$ has {\it constant sectional curvature}
$\kappa$ if ${}^gR=\kappa R_{\id}$.
Let $V$ be a vector space of signature $(r,s)$. For $\varrho>0$, let
$S^\pm(r,s;\varrho)$ be the pseudo-spheres of spacelike and timelike vectors of length
$\pm\varrho^{-1}$:
$$S^\pm(r,s;\varrho):=\{v\in V:(v,v)=\pm\varrho^{-2}\}.$$
The following result is well known so we omit the proof in the interests of brevity.
\begin{lemma}\label{crefb} Let $\varrho>0$. Then $S^\pm(r,s;\varrho)$ has constant sectional
curvature $\pm\varrho$ and is a rank $2$ spacelike, timelike, and mixed Jordan IP
pseudo-Riemannian manifold. Any pseudo-Riemannian manifold of constant
sectional curvature is either flat or is locally isometric to one of these
manifolds.\end{lemma}

It is also possible to construct examples of rank $2$ spacelike Jordan IP pseudo-Riemannian
manifolds by taking twisted products. We introduce the following notational conventions. Let
$\varrho>0$, let $\varepsilon=\pm1$, and let
$\delta=\pm1$ be given. Let
\begin{eqnarray*}
&&M:=I\times S^\delta(r,s;\varrho),\qquad\qquad\quad
 f(t):=\varepsilon\kappa t^2+At+B,\\
&&ds_M^2:=\varepsilon dt^2+f(t)d_{S^\delta(r,s;\varrho)}^2,
 \phantom{aaa} C(t):=f^{-2}\{f\kappa-\textstyle\frac14\varepsilon f_t^2\},\\
&&\phi:=-\id\text{ on }TS^\delta(r,s;\varrho),\qquad\quad
  \phi(\partial_t):=\partial_t.\end{eqnarray*}
Choose $\{\kappa,A,B\}$ so $f\kappa-\frac14\varepsilon
f_t^2\ne0$ or equivalently so that $A^2-4\varepsilon\kappa B\ne0$. Choose the
interval $I$ so $f(t)\ne0$ on $I$.

\begin{lemma}\label{crefc} We have that
$(M,g_M)$ is a rank $2$ spacelike, timelike, and mixed Jordan IP pseudo-Riemannian
manifold with ${}^gR=CR_\phi$.
\end{lemma}

\demo{Proof} Let $m=r+s$.  Fix $P\in S:=S^\delta(r,s;\varrho)$. Choose local coordinates
$x=(x_1,...,x_{m-1})$ for $S$ which are centered at $P$. We let $x_0=t$ to define local
coordinates
$(x_0,...,x_{m-1})$ on $M$. We let indices
$a$,
$b$,
$c$, and
$d$ range from
$1$ to
$m-1$ and index the local coordinate frame for $S$. We let indices $i$, $j$,
$k$, and
$\ell$ range from $0$ to $m-1$ and index the full local coordinate frame for $M$. Let
$g_{ij}$ and
$\tilde g_{ab}$ denote the components of the metric tensors on $M$ and $S$ relative to
this local coordinate frame. We normalize the coordinates on
$S$ so that
$\partial_a\tilde g_{bc}(P)=0$. We have $g_{ab}=f\tilde g_{ab}$, $g_{0a}=0$, and
$g_{00}=\varepsilon$. Let
$f_t:=\partial_tf$ and
$f_{tt}=\partial_t^2f$. Let $\Gamma$ and $\tilde\Gamma$ be the
Christoffel symbols:
\begin{eqnarray*}
&&\Gamma_{ijk}=\textstyle\frac12\{\partial_jg_{ik}+\partial_ig_{jk}-\partial_kg_{ij}\},
\quad\phantom{.}\tilde\Gamma_{abc}=
\textstyle\frac12\{\partial_b\tilde g_{ac}+
    \partial_a\tilde g_{bc}-\partial_c\tilde g_{ab}\},\\
&&\Gamma_{0ab}=\Gamma_{a0b}=\textstyle\frac12f_t\tilde g_{ab},\qquad\qquad\quad
    \phantom{a.}
      \Gamma_{ab0}=-\textstyle\frac12f_t\tilde g_{ab},\quad
\Gamma_{abc}=f\tilde\Gamma_{abc},\\
&&\Gamma_{0a}{}^b=\Gamma_{a0}{}^b\equiv\textstyle\frac12f^{-1}f_t\tilde g_a{}^b
  ,\qquad\qquad
 \Gamma_{ab}{}^0=-\textstyle\frac12\varepsilon f_t\tilde g_{ab},\ 
 \Gamma_{ab}{}^c=\tilde\Gamma_{ab}{}^c;\end{eqnarray*}
the Christoffel symbols where $0$ appears twice or three times vanish. Thus:
\begin{eqnarray*}
   R_{ijk}{}^l&=&\partial_i\Gamma_{jk}{}^l-\partial_j\Gamma_{ik}{}^l
   +\textstyle\sum_n\Gamma_{in}{}^l\Gamma_{jk}{}^n
   -\textstyle\sum_n\Gamma_{jn}{}^l\Gamma_{ik}{}^n,\\
 R_{abc}{}^d&\equiv&\{\partial_a\tilde\Gamma_{bc}{}^d-\partial_b\tilde\Gamma_{ac}{}^d\}
        +\Gamma_{a0}{}^d\Gamma_{bc}{}^0-\Gamma_{b0}{}^d\Gamma_{ac}{}^0\\
    &\equiv&\tilde R_{abc}{}^d
    -\textstyle\frac14\varepsilon f_t^2f^{-1}
      \{\tilde g_a{}^d\tilde g_{bc}-\tilde g_{ac}\tilde g_b{}^d\}
   \text{ mod }O(|x|),\\
 R_{abcd}&\equiv&(f\kappa-\textstyle\frac14\varepsilon f_t^2)
      \{\tilde g_{ad}\tilde g_{bc}-\tilde g_{ac}\tilde g_{bd}\}\\
   &\equiv&f^{-2}(f\kappa-\textstyle\frac14\varepsilon f_t^2)\
      \{g_{ad}g_{bc}-g_{ac}g_{bd}\}\text{ mod }O(|x|).\end{eqnarray*}
Since $S$ is a symmetric space, there is an isometry of $S$ fixing $P$ which acts as $-1$
on $T_PS$. This extends to an isometry of $S$. Thus $R_{0abc}=0$ for any $a$, $b$, and
$c$. Let $a\ne b$. The isotropy group of isometries of $S$ fixing $P$ acts on $T_PM$ as
$O(p,q)$. Thus we can find an isometry of $S$ which sends $e_a\rightarrow e_a$ and
$e_b\rightarrow -e_b$. This implies $R_{0ab0}=0$. We compute the remaining
curvature:
\begin{eqnarray*}
   R_{0aa}{}^0&\equiv&\partial_0\Gamma_{aa}{}^0-\Gamma_{ab}{}^0\Gamma_{0a}{}^b
   \equiv-\textstyle\frac12\varepsilon f_{tt}\tilde g_{aa}
    +\textstyle\frac14\varepsilon f^{-1}f_t^2\tilde g_{ab}\tilde g_a{}^b\\
   R_{0aa0}&\equiv&f^{-2}\varepsilon\{-\textstyle\frac12 ff_{tt}+\frac14f_t^2\}
    \{g_{aa}g_{00}-g_{0a}g_{a0}\}.
\end{eqnarray*}
We wish to find $C(t)$ and $\phi$ so that 
${}^gR=CR_\phi$. If 
$$(f\kappa-\textstyle\frac14\varepsilon f_t^2)=
\varepsilon\{-\textstyle\frac12 ff_{tt}+\frac14f_t^2\},$$
then we may take $\phi=\id$; the resulting metric then has constant sectional curvature.
Since this does not give rise to a new family of metrics, we define instead:
$$\phi(\partial_t)=\partial_t,\text{ and }\phi(\partial_a)=-\partial_a.$$
To ensure that $R=CR_\phi$, we solve the equation:
$$(f\kappa-\textstyle\frac14\varepsilon f_t^2)=-
\varepsilon\{-\textstyle\frac12 ff_{tt}+\frac14f_t^2\}\text{ i.e. }
\kappa=\frac12\varepsilon f_{tt}.$$
This implies that the warping function is quadratic so
$f(t)=\varepsilon \kappa t^2+At+B$. \qed
\enddemo

\section{\label{qctSD}The proof of Theorem \ref{areff}}
Throughout this section, we will let $(M,g)$ be a connected pseudo-Riemannian manifold of
signature
$(p,q)$ with $q\ge5$ and dimension $m:=p+q$ which is rank $2$ spacelike Jordan IP. We use
Theorem
\ref{arefd} to express ${}^gR_P=C(P)R_{\phi(P)}$ where $\phi$ is self-adjoint and $C\ne0$. We
suppose for the moment that ${}^gR$ is never nilpotent so that we can normalize $\phi$ so
$\phi^2=\pm\id$. The maps
$P\rightarrow C(P)$ and $P\rightarrow\phi(P)$ can then be chosen to be smooth, at least
locally. To have a unified notation, we complexify the tangent bundle and extend $\phi$, 
$(\cdot,\cdot)$, and the curvature tensor ${}^gR$ to be complex
multi-linear. Let:
$$\tilde\phi:=\left\{
\begin{array}{ll}
\phantom{\sqrt{-1}}\phi&\text{ if }\phi^2=\phantom{-}\id,\\
\sqrt{-1}\phi&\text{ if }\phi^2=-\id.\end{array}\right.$$
Since $\tilde\phi^2=\id$, the eigenvalues of $\tilde\phi$ are $\pm1$ and $\tilde\phi$ is
diagonalizable. Let
$$\cal{F}^\pm:=\{X\in TM\otimes\mathbb{C}:\tilde\phi X=\pm X\}\text{ and }
  \nu^\pm:=\dim\cal{F}^\pm.$$
If $\nu^-=0$ or $\nu^+=0$, then $\tilde\phi=\phi=\pm\id$ and $(M,g)$ has constant
sectional curvature and assertion (2) of Theorem \ref{areff} holds. Thus, we
suppose
$\nu^-\ge1$ and $\nu^+\ge1$.
Let
$u^\pm\in\cal{F}^\pm$. As
$(u^+,u^-)=(\tilde\phi u^+,u^-)=(u^+,\tilde\phi u^-)=-(u^+,u^-)$,
$\cal{F}^-$ and $\cal{F}^+$ are non-degenerate orthogonal complex distributions. As we
are working over $\mathbb{C}$, we can find a local orthormal frame
$\cal{B}:=\{e_1,...,e_m\}$ for $TM\otimes\mathbb{C}$ so
$$\cal{F}^-=\Pspan\{e_1,...,e_{\nu^-}\},\ \cal{F}^+=\Pspan\{e_{\nu^-+1},...,e_m\},
\text{ and }(e_i,e_j)=\delta_{ij}.$$
Let Roman indices $a$, $b$, etc. range from $1$ to
$\nu^-$ and index the frame for $\cal{F}^-$, let Greek indices
$\alpha$, $\beta$, etc. range from $\nu^-+1$ to $m$ and index the frame for $\cal{F}^+$,
and let Roman indices
$i$,
$j$, etc. range from $1$ to $m$ and index the frame for $TM\otimes\mathbb{C}$. Let
$\tilde\phi_{ij;k}$ be the components of $\nabla\tilde\phi$ relative to such a normalized
basis.
\begin{lemma}\label{drefa} Let $(M,g)$ be a connected pseudo-Riemannian manifold of signature
$(p,q)$ with $q\ge5$ so that ${}^gR_P=C(P)R_{\phi(P)}$ where $\phi(P)^2=\pm\id$ and
$C(P)\ne0$. Assume that $\nu^+\ge1$ and $\nu^-\ge1$. 
\begin{enumerate}
\smallskip\item We have $\tilde\phi_{ij;k}=\tilde\phi_{ji;k}$, $\tilde\phi_{ab;k}=0$,
   $\tilde\phi_{\alpha\beta;k}=0$, and $\tilde\phi_{a\alpha;i}=-2\Gamma_{ia\alpha}$.
\smallskip\item If $i$, $j$, and $k$ are distinct indices, then
   $\tilde\phi_{ij;k}=\tilde\phi_{ik;j}$.
\smallskip\item The only non-zero components of $\nabla\tilde\phi$ are
$\tilde\phi_{a\alpha;a}=\tilde\phi_{\alpha a;a}=-2\Gamma_{aa\alpha}$ and
$\tilde\phi_{a\alpha;\alpha}=\tilde\phi_{\alpha a;\alpha}=-2\Gamma_{\alpha a\alpha}$.
\end{enumerate}
\end{lemma}

\demo{Proof}  We adapt arguments from \cite{refGLS} to prove this result. As the
frame is orthonormal, $\Gamma_{kij}=-\Gamma_{kji}$. Let $\varepsilon_i:=(\phi
e_i,e_i)=\pm1$; $\varepsilon_a=-1$ and $\varepsilon_\alpha=+1$. We prove assertion (1) by
computing:
\begin{eqnarray*}
 \tilde\phi_{ij;k}&=&(\nabla_{e_k}\tilde\phi e_i-\tilde\phi\nabla_{e_k}e_i,e_j)
   =(\nabla_{e_k}\tilde\phi e_i,e_j)-(\nabla_{e_k}e_i,\tilde\phi e_j)\\
    &=&(\varepsilon_i-\varepsilon_j)(\nabla_{e_k}e_i,e_j)
    =(\varepsilon_i-\varepsilon_j)\Gamma_{kij}.\end{eqnarray*}
Let ${}^gR_{ijkl;n}$ be the components of $\nabla{}^gR$. With our normalizations, we can
express:
\begin{eqnarray*}
&&{}^gR_{ijkl}=\varepsilon
  C\{\tilde\phi_{il}\tilde\phi_{jk}-\tilde\phi_{ik}\tilde\phi_{jl}\}\text{ where }
  \phi^2=\varepsilon\id\text{ and},\\
  &&{}^gR_{ijkl;n}=\varepsilon 
C(\tilde\phi_{il;n}\tilde\phi_{jk}+\tilde\phi_{il}\tilde\phi_{jk;n}
         -\tilde\phi_{ik;n}\tilde\phi_{jl}-\tilde\phi_{ik}\tilde\phi_{jl;n})\\
        &&\qquad\qquad\qquad
  +C_{;n}(\tilde\phi_{il}\tilde\phi_{jk}-\tilde\phi_{ik}\tilde\phi_{jl}).\end{eqnarray*}
We use the {\it second Bianchi identity}:
$$0={}^gR_{ijkl;n}+{}^gR_{ijln;k}+{}^gR_{ijnk;l}.$$
Fix distinct indices $i$, $j$, and $k$. As $q\ge5$, we may choose an index
$l$ which is distinct from $i$,
$j$, and $k$. We have $\tilde\phi_{ij}=0$ for
$i\ne j$. Since $C\ne0$, we may prove assertion (2) by computing:
$$0={}^gR_{illj;k}+{}^gR_{iljk;l}+{}^gR_{ilkl;j}
  =\varepsilon C\varepsilon_\ell\tilde\phi_{ij;k}+0
  -\varepsilon C\varepsilon_\ell\tilde\phi_{ik;j}.$$
If $i$, $j$, and $k$ are distinct, then
$\tilde\phi_{ij;k}=\tilde\phi_{ik;j}=\tilde\phi_{jk;i}$. Since at least two of the
indices must index elements of $\cal{F}^-$ or $\cal{F}^+$, $\tilde\phi_{ij;k}=0$
by assertion (1). Thus $(i,j,k)$ is a permutation of $(a,a,\alpha)$ or
$(a,\alpha,\alpha)$. Assertion (3) follows as
$$\tilde\phi_{aa;\alpha}=\tilde\phi_{\alpha\alpha;a}=0\text{ and }
  \tilde\phi_{a\alpha;i}=\tilde\phi_{\alpha a;i}.\ \qed$$\enddemo

We continue our study of $\nabla\tilde\phi$.
\begin{lemma}\label{drefb}Let $(M,g)$ be a connected pseudo-Riemannian manifold of
signature
$(p,q)$ with $q\ge5$ so that ${}^gR_P=C(P)R_{\phi(P)}$ where $\phi(P)^2=\pm\id$ and
$C(P)\ne0$. Assume that $\nu^+\ge1$ and $\nu^-\ge1$.
\begin{enumerate}
\smallskip\item If $\nu^-\ge2$, then $C_{;a}=C\tilde\phi_{\alpha a;\alpha}$
  and $\Gamma_{\alpha a\alpha}=-\frac12\frac{C_{;a}}C$.
\smallskip\item If $\nu^-\ge3$, then $C_{;a}=0$,
   $C_{;\alpha}=-2C\tilde\phi_{a\alpha;a}$, and
    $\Gamma_{aa\alpha}=\frac{C_{;\alpha}}C$.
\smallskip\item If $\nu^+\ge2$, then $C_{;\alpha}=-C\tilde\phi_{a\alpha;a}$ and
    $\Gamma_{aa\alpha}=\frac12\frac{C_{;\alpha}}C$.
\smallskip\item If $\nu^+\ge3$, then $C_{;\alpha}=0$,
   $C_{;a}=+2C\tilde\phi_{a\alpha;\alpha}$, and $\Gamma_{\alpha a\alpha}=-\frac{C_{;a}}C$.
\end{enumerate}\end{lemma}

\demo{Proof} If $\nu^-\ge2$, then we may choose $a\ne b$ and use Lemma \ref{drefa} and the
second Bianchi identity to prove assertion (1) by computing:
$$0={}^gR_{a\alpha\alpha a;b}+{}^gR_{a\alpha ab;\alpha}
     +{}^gR_{a\alpha b\alpha;a}
    =-\varepsilon C_{;b}+\varepsilon C\tilde\phi_{\alpha b;\alpha}+0.
$$
If $\nu^-\ge3$, then we may choose distinct indices $a$, $b$, and
$c$ and use Lemma \ref{drefa} and the second
Bianchi identity to compute
\begin{eqnarray*}
   0&=&{}^gR_{cbbc;a}+{}^gR_{cbca;b}+{}^gR_{cbab;c}
    =\varepsilon C_{;a}+0+0\text{ and}\\
   0&=&{}^gR_{cbbc;\alpha}+{}^gR_{cbc\alpha;b}
      +{}^gR_{cb\alpha b;c}
   =\varepsilon C_{;\alpha}+\varepsilon
C\tilde\phi_{b\alpha;b}
   +\varepsilon C\tilde\phi_{c\alpha;c}.\end{eqnarray*}
Thus $C_{;a}=0$ and
$C_{;\alpha}=-C(\phi_{b\alpha;b}+\phi_{c\alpha;c})$. Similarly
$C_{;\alpha}=-C(\phi_{a\alpha;a}+\phi_{c\alpha;c})$. Thus
$\phi_{b\alpha;b}=\phi_{a\alpha;a}=\phi_{c\alpha;c}$ and assertion (2) follows. 

We use the
same argument to prove assertions (3) and (4) with appropriate changes of sign. 
If $\nu^+\ge2$, then we may choose $\alpha\ne\beta$ to compute:
$$0={}^gR_{\alpha aa\alpha;\beta}+{}^gR_{\alpha a \alpha\beta;a}
     +{}^gR_{\alpha a \beta a;\alpha}
    =-\varepsilon C_{;\beta}-\varepsilon C\tilde\phi_{a\beta;a}+0.
$$
If $\nu^+\ge3$, then we may choose distinct indices $\alpha$, $\beta$, and
$\gamma$ to compute
\begin{eqnarray*}
   0&=&{}^gR_{\gamma\beta\beta\beta\gamma;\alpha}
    +{}^gR_{\gamma\beta\gamma\alpha;\beta}+{}^gR_{\gamma\beta\alpha\beta;\gamma}
    =\varepsilon C_{;\alpha}+0+0\text{ and}\\
   0&=&{}^gR_{\gamma\beta\beta\gamma;a}
     +{}^gR_{\gamma\beta\gamma a;\beta}
      +{}^gR_{\gamma\beta a\beta;\gamma}
   =\varepsilon C_{;a}-\varepsilon
C\tilde\phi_{\beta a;\beta}
   -\varepsilon C\tilde\phi_{\gamma a;\gamma c}.\ \qed\end{eqnarray*}
\enddemo

We can now show that $\phi^2=\id$ if $C\ne0$.
\begin{lemma}\label{drefc} Let $(M,g)$ be a connected pseudo-Riemannian manifold of signature
$(p,q)$ with $q\ge5$ so that ${}^gR_P=C(P)R_{\phi(P)}$ where $\phi(P)^2=\pm\id$ and
$C(P)\ne0$. Assume that $\nu^+\ge1$ and $\nu^-\ge1$.
\begin{enumerate}
\smallskip\item We do not have $\nabla\tilde\phi=0$. Furthermore either $\nu^-\le2$ or
$\nu^+\le2$.
\smallskip\item We have $\phi^2=\id$. Furthermore either $\nu^-=1$ or $\nu^+=1$.
\end{enumerate}\end{lemma}

\demo{Proof} Assume that $\nabla\tilde\phi=0$. We use Lemma \ref{drefa}
to see that
$\Gamma_{ia\alpha}=-\frac12\tilde\phi_{a\alpha;i}=0$.
Thus $\nabla e_a\in\cal{F}^-$ so  ${}^gR(e_a,e_\alpha)e_a\in\cal{F}^-$. This shows:
$$0=^gR(e_a,e_\alpha,e_a,e_\alpha)=\varepsilon C.$$ 
This shows that $C=0$ which is false. Thus $\nabla\tilde\phi\ne0$. Next suppose that
$\nu^-\ge3$ and $\nu^+\ge3$. We may then apply Lemma \ref{drefb} to see that $dC=0$ and
$\nabla\tilde\phi=0$ which is false. This establishes assertion (1).

Suppose that $\phi^2=-\id$. Then $(\phi v,\phi v)=(\phi^2v,v)=-(v,v)$ so
$\phi$ interchanges the roles of spacelike and timelike vectors. Thus
$p=q\ge 5$. Since $\tilde\phi=\sqrt{-1}\phi$, $\tilde\phi$ is purely imaginary. Thus
conjugation interchanges the distributions $\cal{F}^+$ and $\cal{F}^-$ so
$\nu^+=\nu^-=q\ge5$. This contradicts assertion (1). Consequently $\phi^2=\id$.

By replacing $\phi$ by $-\phi$, we may assume
$\nu^-\le\nu^+$. To establish the final assertion, we must rule out the
case
$\nu^-=2$. We may then use assertions (1), (3), and (4) of Lemma \ref{drefb} to see
$dC=\nabla\tilde\phi=0$ which contradicts assertion (1).
\qed\enddemo

Assertion (1) of Theorem \ref{areff} follows from Lemma \ref{drefc} and from:

\begin{lemma}\label{drefd} Let $(M,g)$ be a connected pseudo-Riemannian manifold of
signature $(p,q)$ for $q\ge5$. Assume that $(M,g)$ is rank $2$ spacelike Jordan IP.
Then either ${}^gR$ is nilpotent for all points of $M$ or ${}^gR$ is nilpotent at no point
of $M$.
\end{lemma}

\demo{Proof} Let $S_1$ be the set of all points of $M$ where ${}^gR$ is
nilpotent and let $S_2$ be the complementary subset of all points of $M$ where ${}^gR$
is not nilpotent. We assume that both $S_1$ and $S_2$ are non-empty and argue for a
contradiction. Since $S_1$ is closed, since $S_2$ is open, and since $M$ is connected, we
may conclude that $S_2$ is not closed. Thus we may choose points $P_i\in S_2$ so that
$P_i\rightarrow P_\infty\in S_1$. Let $R_n$ be the
curvature tensor at $P_n$ and let $V_n$ be the tangent space at $P_n$.

Let $\cal{B}:=\{f_1^-,...,f_p^-,f_1^+,...,f_q^+\}$ be a
normalized local orthonormal frame for the tangent bundle near $P_\infty$ and let $V^\pm$
be the associated maximal spacelike and timelike distributions. We set
$V_n^\pm:=V^\pm(P_n)$.  Since $P_n\in S_2$, we can express
$R_n=C_nR_{\phi_n}$. Since $\text{Tr}\{R_n(\pi)^2\}=-C_n^2$ for any spacelike $2$ plane
$\pi$, we have $C_n\rightarrow 0$ since $R_\infty$ is nilpotent.

 We apply Lemma
\ref{drefc} to see
$\phi_n^2=\id$ and thus we do not need to complexify to define the sub-bundles
$\cal{F}^\pm_n$ of $V_n$. By replacing
$\phi_n$ by
$-\phi_n$ if necessary, we may apply Lemma \ref{drefc} to see $\dim\cal{F}^-_n\le1$. Consequently
$\dim\cal{F}^+\ge p+q-1$ so $\dim\cal{F}^+\cap V_n^+\ge q-1\ge 4$. Choose elements
$v_n,w_n\in V_n^+\cap\cal{F}_n^+$ so that
$\{v_n,w_n\}$ forms an orthonormal spacelike set. We choose a compact neighborhood of
$P_\infty$ over which
$S(V^+)$ is compact. By passing to a subsequence,
we can suppose that
$v_n\rightarrow v_\infty$ and $w_n\rightarrow w_\infty$. Let $\pi_n:=\Pspan\{v_n,w_n)$.
Then we have that $\pi_n\rightarrow\pi_\infty:=\Pspan\{v_\infty,w_\infty\}$. The planes $\pi_n$
and
$\pi_\infty$ are spacelike. Let
$z_\infty\in V_\infty$. Choose elements $z_n\in V_n$ so $z_n\rightarrow z_\infty$. As
$\phi_nv_n=v_n$ and
$\phi_nw_n=w_n$, we have:
\begin{eqnarray*}
R(\pi_\infty)z_\infty
 &=&\lim_{n\rightarrow\infty}R_n(v_n,w_n)z_n\\
&=&\lim_{n\rightarrow\infty}C_n\{(\phi w_n,z_n)\phi v_n-
(\phi v_n,z_n)\phi w_n\}\\
&=&\lim_{n\rightarrow\infty}C_n\{(w_n,z_n)v_n-(v_n,z_n)w_n\}\\
  &=&\{\lim_{n\rightarrow\infty}C_n\}\cdot
   \{(w_\infty,z_\infty)v_\infty-(v_\infty,z_\infty)w_\infty\}=0.\end{eqnarray*}
This contradicts the assumption that $R$ has spacelike rank $2$.
\qed\enddemo

Since $\phi^2=\id$ so the distributions $\cal{F}^\pm$ are real. By replacing
$\phi$ by $-\phi$ if necessary, we may suppose that $\nu^-\le\nu^+$ and thus $\nu^-\le1$.
If $\nu^-=0$, then $(M,g)$ has constant sectional curvature and assertion (2) of Lemma
\ref{areff} holds. We therefore suppose $\nu^-=1$. The distributions $\cal{F}^\pm$ are
non-degenerate. Let
$\{e_1\}$ be a local orthonormal section to $\cal{F}^-$ and let $\{e_2,...,e_m\}$ be a
local orthonormal frame for $\cal{F}^+$; since we are working over $\mathbb{R}$, we no longer
impose the normalization that
$(e_i,e_i)=+1$. By replacing $g$ by $-g$, we may assume without loss of generality that
$e_1$ is spacelike. If $\alpha\ne\beta$, then we may apply Lemma \ref{drefa} to compute:
\begin{eqnarray*}
([e_\alpha,e_\beta],e_1)&=&\Gamma_{\alpha\beta1}-\Gamma_{\beta\alpha1}
    =\Gamma_{\beta1\alpha}-\Gamma_{\alpha1\beta}\\
    &=&-\textstyle\frac12(\tilde\phi_{1\alpha;\beta}-\tilde\phi_{1\beta;\alpha})=0.
\end{eqnarray*}
Thus the foliation $\cal{F}^+$ is integrable. 
Let
$y=(y^1,...,y^{m-1})$ be local coordinates on a leaf of this foliation. 
We define geodesic tubular coordinates on $M$ by setting:
$$T(t,y):=\operatorname{exp}_y(te_1(y)).$$

Assertion (3) of Theorem \ref{areff} will follow from:

\begin{lemma}\label{drefe}  Let $(M,g)$ be a connected pseudo-Riemannian manifold of
signature $(p,q)$ for $q\ge5$. Assume that ${}^gR=CR_\phi$ where $\phi^2=\id$ and $C\ne0$.
Assume that $\nu^-=1$ and that $\cal{F}^-$ is spacelike.
\begin{enumerate}
\smallbreak\item For fixed $y_0$, the curves $t\rightarrow T(t,y_0)$ are unit
speed geodesics in $(M,g)$ which are perpendicular to the leaves of the foliation
$\cal{F}^+$.
\smallbreak\item For fixed $t_0$, the surfaces $T(t_0,y)$ are leaves of the
foliation $\cal{F}^+$ and inherit metrics of constant sectional curvature.
\smallbreak\item Locally $ds^2=dt^2+fds_\kappa ^2$ where $f(t)$ is non-zero smooth
function and $ds_\kappa^2$ is a metric of constant sectional curvature $\kappa$.
\smallbreak\item The warping function $f(t)=\kappa t^2+At+B$ where $A^2-4\kappa B\ne0$.
\end{enumerate}\end{lemma}

\demo{Proof} We choose a local orthonormal frame $\{e_i\}$ for $TM$ so $e_1$ spans
$\cal{F}^-$ and $\{e_2,...,e_m\}$ spans $\cal{F}^+$. We set
$\varepsilon_i:=(e_i,e_i)=\pm1$; by assumption $\varepsilon_1=+1$. We have
$\dim(\cal{F}^+)=m-1\ge3$. Taking into account the fact that $(e_\alpha,e_\alpha)=\pm1$, we
may apply Lemmas
\ref{drefa} and \ref{drefb} to see that
$$C_{;\alpha}=0,\ 
   \Gamma_{\alpha1\beta}=-\textstyle(e_\alpha,e_\beta)\frac{C_{;1}}C,\text{ and }
   \Gamma_{11\alpha}=\frac12\frac{C_{;\alpha}}C=0.$$
Let $\gamma(t,y_0)$ be an integral curve for $e_1$ starting at a point $y_0$ on the
leaf of the foliation $\cal{F}^+$. Since $e_1$ is a unit
vector, $\Gamma_{111}=0$. As $\Gamma_{11\alpha}=0$, $\gamma$ is a geodesic.
Thus $\gamma(t,y_0)=T(t,y_0)$ so $\partial_t=e_1$. We
compute
$$\partial_t(\partial_t,\partial_\alpha^y)
  =(\partial_t,\nabla_{\partial_t}\partial_\alpha^y)
  =(\partial_t,\nabla_{\partial_\alpha^y}\partial_t)
  =\textstyle\frac12\partial_\alpha^y(\partial_t,\partial_t)=0.$$
This shows $(\partial_t,\partial_\alpha^y)=0$ so
the $\partial_\alpha^y$ span the perpendicular distribution $\cal{F}^+$ and the manifolds
$T(t_0,y)$ are leaves of the foliation $\cal{F}^+$. Since
$C_{;\alpha}=0$, $C$ is constant on
the leaves of $\cal{F}^+$. Let
$\tilde R$ be the curvature of the induced metric on the leaves of the
foliation $\cal{F}^+$. We show that
$\tilde R$ has constant sectional curvature by computing:
\begin{eqnarray*}
  \tilde R(e_\alpha,e_\beta,e_\gamma,e_\sigma)&=&R(e_\alpha,e_\beta,e_\gamma,e_\sigma)
 +\Gamma_{\alpha1\sigma}\Gamma_{\beta\gamma}{}^1
 -\Gamma_{\beta1\sigma}\Gamma_{\alpha\gamma}{}^1\\
 &=&\{C(t)+\textstyle(\frac{C_{;1}}C)^2\}
   \{(e_\beta,e_\gamma)(e_\alpha,e_\sigma)-
    (e_\alpha,e_\gamma)(e_\beta,e_\sigma)\}.\end{eqnarray*}
Let
$\partial_\alpha^y=\Sigma_\gamma a_{\alpha\gamma}e_\gamma$. 
We show that the metric is a warped product by computing:
\begin{eqnarray*}
 (\nabla_{\partial_t}\partial_\alpha^y,\partial_\beta^y)&=&
     (\nabla_{\partial_\alpha^y}\partial_t,\partial_\beta^y)
    =\Sigma_{\gamma\sigma}a_{\alpha\gamma}a_{\beta\sigma}
   (\nabla_{e_\gamma}\partial_t,e_\sigma)\\
 &=&\Sigma_{\gamma\sigma}a_{\alpha\gamma}a_{\beta\sigma}
    (e_\gamma,e_\sigma)\textstyle\frac{C_{;1}}C=
   \textstyle\frac{C_{;1}}Cg_{\alpha\beta}\text{ so}\\
 \partial_tg_{\alpha\beta}&=&(\nabla_{\partial_t}\partial_\alpha^y,\partial_\beta^y)
      +(\partial_\alpha^y,\nabla_{\partial_t}\partial_\beta^y)=\textstyle
    \frac{2C_{;1}}Cg_{\alpha\beta}.\end{eqnarray*}
The
argument used to prove Lemma
\ref{crefc} now shows the warping function is quadratic.
\qed\enddemo

\section{\label{qctSE}Nilpotent spacelike Jordan IP pseudo-Riemannian
manifolds}
The rank $2$ spacelike Jordan IP pseudo-Riemannian manifolds whose curvature operators
are not nilpotent for at least one point of $M$ are
classified in Theorem
\ref{areff}. In this section, we study the remaining case and present some preliminary results.
We focus our attention on the balanced setting $p=q$.

\begin{lemma}\label{erefa} Let $(M,g)$ be a connected 
pseudo\-Riemannian manifold of signature $(p,p)$ for $p\ge5$. Assume that $(M,g)$
is spacelike rank 2 nilpotent Jordan IP.
\begin{enumerate}
\smallskip\item We have ${}^gR=\pm R_\phi$ where $\phi$ is self-adjoint and where
$\ker\phi=\range\phi$.
\smallskip\item $\range\phi$ is an integrable distribution of the
tangent bundle of $M$.
\end{enumerate}\end{lemma}

\demo{Proof} We apply Theorem \ref{arefd} to write $R=\pm R_\phi$ where we normalize
$\phi$ by requiring that $C=\pm1$. The map $P\rightarrow\phi$ can then be chosen to vary
smoothly with
$P$, at least locally. Since $\phi$ contains no spacelike vectors, $\dim\{\ker(\phi)\}\le
p$. Since $\phi$ is self-adjoint and $\phi^2=0$, we show
$\range(\phi)=\ker(\phi)$ and complete the proof of the first assertion by computing:
\begin{eqnarray*}
&&\range(\phi)\subset\ker(\phi),\\
&&\dim\{\range(\phi)\}\le\dim\{\ker(\phi)\},\text{ and}\\
&&2p=\dim\{\range(\phi)\}+\dim\{\ker(\phi)\}\le2\dim\{\ker(\phi)\}\le2p.\end{eqnarray*}

Let $\cal{K}:=\range\phi$ and let
${\cal{S}}$ be a maximal local spacelike distribution. Since ${\cal{K}}$ is totally isotropic,
${\cal{K}}\cap{\cal{S}}=\{0\}$. Thus ${\cal{K}}=\range\phi=\phi{\cal{S}}$ and
$TM={\cal{S}}\oplus\phi{\cal{S}}$.
Let $L(\cdot,\cdot):=(\phi\cdot,\cdot)$ be the associated bilinear form. We have:
$$(\phi s,\phi\tilde s)=(s,\phi^2\tilde s)=0.$$ 
Thus if $0\ne s\in{\cal{S}}$, then there must
exist $\tilde s\in{\cal{S}}$ so that $(\phi s,\tilde s)\ne0$. Thus $L$ is a
non-degenerate bilinear form. We complexify and choose a frame $\{s_a\}$ for
${\cal{S}}$ so:
$$(\phi e_a,e_b)=\delta_{ab}.$$
Let Roman indices $a$, $b$, etc. range from $1$ to $p$ and index this frame for
${\cal{S}}$. Let Greek indices $\alpha$, $\beta$, etc. range from $p+1$ to $2p$ and
index the frame $e_\alpha:=\phi(e_{\alpha-p})$ for ${\cal{K}}$. Let Roman indices $i$, $j$,
etc. range from $1$ to $2p$ and index the frame $\{e_1,...,e_p,\phi
e_1,...,\phi e_p\}$. We have $\phi_{ab}=\delta_{ab}$, $\phi_{a\beta}=0$, and
$\phi_{\alpha\beta}=0$. Let $\alpha$, $\beta$, and $i$ be indices which need not be
distinct. Choose $a\ne i$ and use the second Bianchi identity to compute:
\begin{eqnarray*}
R_{ijkl;n}&=&\phi_{il;n}\phi_{jk}+\phi_{il}\phi_{jk;n}-\phi_{ik;n}\phi_{jl}
   -\phi_{ik}\phi_{jl;n},\\
   0&=&R_{a\alpha\beta a;i}+R_{a\alpha i\beta;a}
    +R_{a\alpha ai;\beta}
   =\phi_{\alpha\beta;i}+0-\phi_{\alpha i;\beta},\text{ and}\\
   \phi_{\alpha i;\beta}&=&\phi_{\alpha\beta;i}
     =(\nabla_{e_i}\phi(\phi e_a),\phi e_b)
   -(\phi\nabla_{e_i}\phi e_a,\phi e_b)=0.\end{eqnarray*}
We clear the previous notation. Let $\alpha=a+p$, $\beta=b+p$, and $\gamma=c+p$ be
indices which are not necessarily distinct. We show that $\cal{K}$ is an integrable
distribution by showing that:
\begin{eqnarray*}
  &&([e_\alpha,e_\beta],e_\gamma)=
   (\nabla_{e_\alpha}\phi e_b,\phi e_c)-(\nabla_{e_\beta}\phi e_a,\phi e_c)\\
   &=&((\nabla_{e_\alpha}\phi-\phi\nabla_{e_\alpha})e_b,\phi e_c)
     -((\nabla_{e_\beta}\phi-\phi\nabla_{e_\alpha})e_a,\phi e_c)\\
   &=&\phi_{b\gamma;\alpha}-\phi_{a\gamma;\beta}=0\text{ so }
  [e_\alpha,e_\beta]\in C^\infty({\cal{K}}^\perp)=C^\infty({\cal{K}}).\ \qed\end{eqnarray*}
\enddemo

We now construct a nilpotent rank $2$ Jordan IP
pseudo-Riemannian manifold. 

\begin{lemma}\label{erefb} Let $\{e_1,...,e_p,\tilde e_1,...,\tilde e_p\}$ be
a basis for a vector space $V$ and let $\{x^1,...,x^p,\tilde x^1,...,\tilde x^p\}$
be the corresponding dual basis for $V^*$. Define an inner product of signature
$(p,p)$ on $V$ by $(e_i,\tilde e_j)=\delta_{ij}$ and
$(e_i,e_j)=(\tilde e_i,\tilde e_j)=0$ for all $i,j$. Give $W:=V\oplus\mathbb{R}$ the direct
sum inner product which has signature
$(p,p+1)$. Let
$f(\tilde x_1,...,\tilde x_p)$ be a real valued function so that $df(0)=0$ and so that
$\det(\partial_i^{\tilde x}\partial_j^{\tilde x}f)(0)\ne0$. The embedding $F:=\id\oplus f$
of $V$ into
$W$ defines a hypersurface $(M,g)$ of $W$ so that $(M,g)$ is a spacelike and timelike rank
$2$ Jordan IP nilpotent pseudo-Riemannian manifold of signature $(p,p)$ close to the
origin.
\end{lemma}

\demo{Proof} Since $df(0)=0$, the natural identification of $T_0M$ with $V$ is an
isometry. Thus the metric is non-degenerate and has signature $(p,p)$ sufficiently close
to the origin. Let $L_P$ and $S_P$ be the associated second fundamental form
and shape operators at $P$. Then Lemma \ref{crefa} shows that ${}^gR=R_{S_P}$.
Since
$\partial_i^xF=0$,
$L_P(\partial_i^x,*)=0$ for any point of the manifold. Thus $S_P(\partial_i^x)=0$ so
$$\Pspan\{\partial_i^x\}\subset\ker(S_P).$$
Since the Hessian of $f$ is non-singular at $0$, $S_0$ is invertible on
$\Pspan\{\partial_i^{\tilde x}\}$. Thus $\dim\ker(S_0)=p$ and hence by shrinking the
neighborhood $\cal{O}$ if necessary we may suppose that
$\dim\ker(S_P)\le p$ for $P\in\cal{O}$. It now follows that
$$\ker(S_P)=\Pspan\{\partial_i^x\}\text{ for all }P\in\cal{O}.$$
Since $F_*(\partial_i^x)=e_i$, we have $g(\partial_i^x,\partial_j^x)=0$ for all $i,j$ and
hence $\ker(S_P)$ is totally isotropic and in particular $\ker(S_P)$ contains no spacelike
or timelike vectors. Since $\dim\{\ker(S_P)\}=p$, $\ker(S_P)=\ker(S_P)^\perp=\range(S_P)$
so $S_P^2=0$. Thus ${}^gR$ is nilpotent.\qed\enddemo

\end{document}